\documentclass{article}
\usepackage{mathptmx}      
\usepackage{amsthm,amsfonts,amsmath}
\usepackage{amscd,amssymb,latexsym,epsfig}
\usepackage[all]{xy}
\title{The backward $\lambda$-Lemma and Morse filtrations
       \\\vspace{.5cm}
       {\small
       Proceedings
       \\
       {\it IX Nonlinear differential equations}
       \\
       on the occasion of the 60th birthday of Bernhard Ruf
       \\
      \vspace{.2cm} 
      17-21 September 2012
       \\
       Jo\~{a}o Pessoa, Para\'{\i}ba, Brazil
       \\
       }
       }
\date{\vspace{-5ex}}
\author{Joa Weber\footnote{
                 Financial support:
                 FAPESP grants 2011/01830-1 and 2012/11634-8.
                 Author address:
                 IMECC,
                 Rua S\'{e}rgio Buarque de Holanda 651,
                13083-859 Campinas, SP, Brasil.
                 joa@math.sunysb.edu.
                 MSC
                 58-02 (Primary)
                 58B05 35K91 (Secondary).
                 Date: 16 Mar 2013. 
                 }
                 \\
        UNICAMP\\
       }
\newtheorem{theorem}{Theorem}[section]

\theoremstyle{definition}
\newtheorem{definition}[theorem]{Definition}

\newtheorem{remark}[theorem]{Remark}

\theoremstyle{remark}

\newcommand{\N}{{\mathbb{N}}}

\newcommand{\R}{{\mathbb{R}}}

\newcommand{\Z}{{\mathbb{Z}}}
\newcommand{\Bb}{{\mathcal{B}}}
\newcommand{\Dd}{{\mathcal{D}}}

\newcommand{\Ff}{{\mathcal{F}}}
\newcommand{\Gg}{{\mathcal{G}}}   

\newcommand{\Ss}{{\mathcal{S}}}

\newcommand{\Ww}{{\mathcal{W}}}

\newcommand{\dist}{{\rm dist}}     
\newcommand{\IND}{{\rm ind}}       

\newcommand{\grad}{{\rm grad }}    
\newcommand{\Crit}{{\rm Crit}}        
\newcommand{\Ho}{{\rm H}}             
\newcommand{\CM}{{\rm CM}}            
\newcommand{\HM}{{\rm HM}}            

\newcommand{\D}{{\rm D}}

\newcommand{\eps}{{\varepsilon}}

\def\NABLA#1{{\mathop{\nabla\kern-.5ex\lower1ex\hbox{$#1$}}}}
\def\Nabla#1{\nabla\kern-.5ex{}_{#1}}
\def\Tabla#1{\Tilde\nabla\kern-.5ex{}_{#1}}

\def\Abs#1{\left|#1\right|}            

\def\Norm#1{\left\|#1\right\|}

\renewcommand{\Tilde}{\widetilde}

\newcommand{\p}{{\partial}}

\begin{document}

\maketitle


\begin{abstract}
Consider the infinite dimensional hyperbolic
dynamical system provided by the (forward) heat semi-flow
on the loop space of a closed Riemannian manifold $M$.
We use the recently discovered backward $\lambda$-Lemma
and elements of Conley theory to construct
a Morse filtration of the loop space whose cellular
filtration complex represents the Morse complex
associated to the downward $L^2$-gradient of the classical
action functional.
This paper is a survey. Details and proofs will
be given in~\cite{Joa-CONLEY}.
\end{abstract}

\tableofcontents

\section{Introduction}

Consider a closed smooth manifold $M$ of dimension
$n\ge 1$ equipped with a Riemannian metric
and the Levi-Civita connection $\nabla$.
Pick a smooth function $V:S^1\times M$ and set
$V_t(q):=V(t,q)$. Here and throughout we identify
$S^1=\R/\Z$.

For smooth maps
$\R\times S^1\to M:(s,t)\mapsto u(s,t)$
consider the \emph{heat equation}
\begin{equation}\label{eq:heat}
   \p_su - \Nabla{t}\p_tu - \nabla V_t(u) = 0.
\end{equation}
It corresponds to the downward $L^2$-gradient equation
of the \emph{action} given by
$$
     \Ss_V(x) = \int_0^1 \left(
     \frac12\Abs{\dot x(t)}^2
     -V(t,x(t)) \right)
     dt
$$
for any element $x:S^1\to M$ of the \emph{free loop space}
$\Lambda M:=W^{1,2}(S^1,M)$ consisting
of absolutely continuous loops in $M$.
The critical points of $\Ss_V$
are the solutions $x\in\Lambda M$ of the ODE
$-\Nabla{t}\dot x-\nabla V_t(x)=0$,
that is the (perturbed) closed geodesics.
\emph{Throughout this paper} we fix a regular value $a$ of
$\Ss_V$ and assume that the Morse-Smale condition holds true
below level $a$. Consider the sublevel set
$\Lambda^a M:=\{\Ss_V<a\}$. In this case the
action is a Morse function on $\Lambda^a M$
and the set of solutions to~(\ref{eq:heat})
that converge to critical points
$x^\pm\in\Lambda^a M$, as $s\to\pm\infty$,
carries the structure of a smooth manifold
whose dimension is given by the Morse index
difference $\IND_V(x)-\IND_V(y)$.
Moreover, the number $n_a$ of elements of the
set $\Crit$ of critical points of $\Ss_V$ in
$\Lambda^a M$ is finite.
By $\Crit_k$ we denote the set
of critical points in $\Lambda^a M$
of Morse index $k$. For each
$x\in \Crit$ pick an orientation of the largest subspace
$E_x$ of the Hilbert space 
$$
     X:=T_x\Lambda M=W^{1,2}(S^1,x^*TM)
$$
on which the Hessian of $\Ss_V$ at $x$ is negative definite.
(The dimension of $E_x$ is finite and called the
\emph{Morse index} of $x$.)

\subsubsection*{Heat flow homology~\cite{Joa-HEATMORSE-II}}

By definition the \emph{Morse chain groups}
$\CM_k=\CM_k(\Lambda^a M,\Ss_V;\Z)$
are the free abelian groups generated by
the (perturbed) closed geodesics $x$ of Morse index $k$
and below level $a$, that is $\Z^{\Crit_k}$.
Set $\CM_k=\{0\}$ in case of the
empty set. The chosen orientations provide
the \emph{characteristic sign} $n_u\in\{\pm 1\}$
for each heat flow solution $u$ of~(\ref{eq:heat})
between critical points of index difference one.
Up to shift in the time variable $s$,
there are only finitely many such $u$.
Counting them with signs $n_u$ provides the
\emph{Morse boundary operator} $\p_k:\CM_k\to\CM_{k-1}$.
By $\HM_k$ we denote the corresponding homology groups.

\subsubsection*{Main result: The natural isomorphism to singular homology~\cite{Joa-CONLEY}}
The idea to use cellular filtrations to calculate Morse
homology goes back at least to Milnor~\cite{MILNOR-h}.
One needs to construct a cellular filtration
$\Ff$ of $\Lambda^a M$ whose cellular filtration
complex $({\rm C}_*\Ff,\p_*)$ precisely represents the
Morse complex, up to natural identification.
In this case we are done, since
\begin{equation}\label{eq:natural-isomorphism}
     \HM_k
     \equiv\Ho_*\left(({\rm C}_*\Ff,\p_*)\right)
     \simeq\Ho_*(\Lambda^a M)
\end{equation}
where the isomorphism is
provided by algebraic topology
given any cellular filtration of $\Lambda^a M$
(related to the Morse complex or not);
see e.g.~\cite{DOLD}.

\section{Morse filtrations and Conley pairs}

\begin{definition}[Cellular filtration and homology]
Assume $\Ff=(F_{-1}\subset F_0\subset F_1\subset\dots
\subset F_\mu)$ is a nested sequence of open subsets
of $\Lambda^a M$ such that
relative singular homology
$H_\ell(F_k,F_{k-1})$ is trivial whenever
$\ell\not= k$ and where $F_{-1}:=\emptyset$.
In this case $\Ff$ is
a \emph{cellular filtration of $\Lambda^a M$}.
For the algebraic topology used in this section
we refer to~\cite{DOLD}. The
\emph{cellular chain complex}
consists of the \emph{cellular chain groups}
${\rm C}_k\Ff:=\Ho_k(F_k,F_{k-1})$
together with the triple boundary operators
$\p_k:\Ho_k(F_k,F_{k-1})\to\Ho_{k-1}(F_{k-1},F_{k-2})$.
A cellular filtration $\Ff$ is called a
\emph{Morse filtration}, if
${\rm C}_k\Ff=\CM_k$ for every $k\in\N$, that is
each relative homology group $H_k(F_k,F_{k-1})$
is generated precisely
by the critical points of Morse index $k$.
\end{definition}
\begin{remark}\label{rem:Morse-filtration}
To establish~(\ref{eq:natural-isomorphism}) we need to
a)~construct
a Morse filtration $\Ff$ of $\Lambda^a M$ and
b)~show that the
associated triple boundary operator counts
heat flow lines according to their characteristic signs
between critical points of index difference one.
How to solve these two problems
is known for flows; cf.~\cite{MILNOR-h}
or~\cite[thm.~2.11]{AM-LecsMCInfDimMfs}.
The solution to b)~carries
over to our semi-flow situation,
since restricted to the (finite dimensional) unstable manifolds
the semi-flow turns into a flow.
It remains to
construct a Morse filtration $\Ff$ of $\Lambda^a M$.
\end{remark}

\subsubsection*{The Abbondandolo-Majer construction for flows~\cite{AM-LecsMCInfDimMfs}}
In their construction of a Morse filtration $\Ff^\prime$ of
$\Lambda^a M$ openness of the sets $F_k^\prime$ follows from
openness of the time-T-map and the Morse property
is a consequence of forward flow invariance of the
open sets $F_k^\prime$. Start by setting $N_0$ equal
to the union of open \emph{local} sublevel sets,
one for each local minimum $x_0$.
Set $F_0^\prime:=N_0$.
Next choose a small
open ball about each index one critical
point and denote their (disjoint) union by
$N_1^\prime$.
Then take the union of $F_0^\prime$ and the whole
forward flow of $N_1^\prime$ and call it
$F_1^\prime:=F_0^\prime\cup\varphi_{[0,\infty)}N_1^\prime$.
Similarly define $F_2^\prime$ and
$F_3^\prime,\ldots F_{n_a}^\prime$.

\subsubsection*{A construction for semi-flows using Conley pairs~\cite{Joa-CONLEY}}
The Cauchy problem associated to the heat
equation~(\ref{eq:heat}) for maps
$[0,\infty)\to \Lambda^a M:s\mapsto u_s=u(s,\cdot)$ is well
posed and leads to the continuous \emph{semi-flow}
\begin{equation*}
     \varphi:[0,\infty)\times\Lambda^a M\to\Lambda^a M
\end{equation*}
called the \emph{heat flow}. In fact $\varphi$
is of class $C^1$ on $(0,\infty)$.
A characteristic feature of the heat flow
is its extremely regularizing nature, namely
$\varphi_s\gamma\in C^\infty(S^1,M)$ whenever
$\gamma \in\Lambda M$ and $s>0$.
Observe that the set of nonsmooth elements is
dense\footnote{
     Pick $\gamma\in\Lambda M$ and a nonsmooth
     $\xi\in W^{1,2}(S^1,x^*TM)$. For large integers $j$
     set $\exp_\gamma (\frac{1}{j}\xi)$.
   }
in $\Lambda M$.
Hence $\varphi_s$ is not an open map for $s>0$ and the
Abbondandolo-Majer method does not work.
Instead we propose the following construction.

It is a very simple---but far reaching---observation that
\emph{by continuity of $\varphi_s$ preimages of open sets are open}.
Define $N_0$ as above.
Observe that the preimage $(\varphi_T)^{-1}N_0$
is open and semi-flow invariant.
Pick any index one critical point $x_1$.
The (one-dimensional) unstable manifold of $x_1$
necessarily\footnote{
      By Palais-Smale and $\Ss_V$ being Morse
      $\gamma_\infty:=\lim_{s\to\infty}\varphi_s\gamma$
      always exists and lies in $\Crit$.
      If $\gamma\in W^u(x_1)$ and $\gamma\not= x_1$, then
      $\gamma_\infty\in\Crit_0$ by Morse-Smale.
    }
enters $N_0$. Consequently our preimage gets very close to
$x_1$ for $T$ very large, however, it never contains $x_1$.
To get over the barrier $x_1$ assume we had an open
neighborhood $N_{x_1}$ of $x_1$ containing no other
critical points and an open subset
$L_{x_1}\subset N_{x_1}$ whose closure does not
contain $x_1$. Assume further that
$L_{x_1}$ is semi-flow invariant in $N_{x_1}$
and every element leaving $N_{x_1}$ under the semi-flow
necessarily runs through $L_{x_1}$ first.
Such a pair $(N_{x},L_{x})$ is called a \emph{Conley pair}
for $x\in\Crit$ and $L_{x}$ is called an \emph{exit set}
for the \emph{Conley set} $N_{x}$.

Pick $x\in\Crit$ and set $c:=\Ss_V(x)$.
For $\eps>0$ small and $\tau>0$ large the sets
\begin{equation}\label{eq:Conley-pair}
\begin{split}
     N_x=N_x^{\eps,\tau}
    &:=\left\{\gamma\in\Lambda^{c+\eps}M\mid
     \Ss_V(\varphi_\tau\gamma)>c-\eps\right\}_x
    \\
     L_x=\, L_x^{\eps,\tau}
    &:=\left\{\gamma\in N_x\mid
     \Ss_V(\varphi_{2\tau}\gamma)<c-\eps\right\}
\end{split}
\end{equation}
form a Conley pair for $x$. Here $\{\ldots\}_x$ indicates
the path connected component that contains $x$.
By Theorem~\ref{thm:inv-fol}~(d) the sets
$N_x$ corresponding to different critical points $x$
are pairwise disjoint.
\begin{figure}
  \centering
  \includegraphics{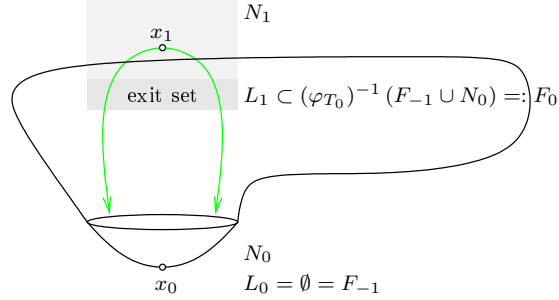}
  \caption{
           Morse filtration $\Ff=\left(\emptyset\subset
           F_0\subset
           F_1\subset\dots\subset F_{n_a}=\Lambda^a M\right)$
          }
  \label{fig:fig-Morse-filtration-JP}
\end{figure}
For $k\in\{0,\dots,n_a\}$ set
$$
     N_k:=\cup_{x\in\Crit_k} N_x
     ,\qquad
     L_k:=\cup_{x\in\Crit_k} L_x.
$$
Consider the sets
\begin{equation}\label{eq:F_k}
     F_{k}
     :=(\varphi_{T_k})^{-1}\left(F_{k-1}\cup N_k\right)
     \supset L_{k+1},\qquad
     k=0,\ldots,n_a-1,
\end{equation}
where the constant $T_k$ is chosen sufficiently
large\footnote{
    Here Palais-Smale, Morse-Smale
    on neighborhoods, and $\Ss_V$ being bounded below enter.
  }
such that the inclusion holds true;
see Figure~\ref{fig:fig-Morse-filtration-JP}.
Because there are no critical points in the complement of
$F_{n_a-1}\cup N_{n_a}$ in $\Lambda^a M$, there is a constant
$T_{n_a}$ such that $\Lambda^a M$ is equal to
$
     F_{n_a}
     :=(\varphi_{T_{n_a}})^{-1}\left(N_{n_a}\cup F_{n_a-1}\right)
$.
Observe that each set $F_k$ is open, because $N_k$ and
$F_{k-1}$ are. Furthermore, although $N_k$ is \emph{not}
semi-flow invariant the union $N_k\cup F_{k-1}$ \emph{is},
because the exit set $L_k$ of $N_k$ is contained in $F_{k-1}$.
Openness and semi-flow invariance heavily enter the
calculation~(\ref{eq:calc}) in the proof of
the Morse filtration property.

\subsubsection*{Morse filtration property}
Constructing suitable homotopy equivalences and
applying the excision axiom of relative homology
one shows that
\begin{equation}\label{eq:calc}
     \Ho_\ell(F_k,F_{k-1})
     \simeq
     \Ho_\ell(N_k,L_k)
     \simeq
     \bigoplus_{x\in\Crit_k}
     \Ho_\ell(N_x,L_x).
\end{equation}
Here the final step uses that $N_k$ is a union of
pairwise disjoint sets $N_x$.
So in order to prove that the nested
sequence $\Ff$ consisting of the open
semi-flow invariant
sets $F_k$ defined by~(\ref{eq:F_k}) is a Morse
filtration of $\Lambda^a M$---thereby concluding
the proof of~(\ref{eq:natural-isomorphism})
via Remark~\ref{rem:Morse-filtration}---it
remains to show that
\begin{equation}\label{eq:claim}
     \Ho_\ell(N_x,L_x)
     \simeq\Ho_\ell(\D^k,\p\D^k)
     \simeq
     \begin{cases}
        \Z& ,\ell=k,
        \\
        0&  \text{, otherwise,}
     \end{cases}
\end{equation}
for every $x\in\Crit_k$.
To prove the first isomorphism was precisely the problem
which inspired us to come up with the backward
$\lambda$-Lemma in~\cite{Joa-LAMBDA}:
Since the part of $N_x$ in the unstable manifold $W^u(x)$ is
a $k$-disk and the corresponding
part of $L_x$ is homotopy equivalent to the disk boundary,
it remains to deformation retract
$(N_x,L_x)$ to its part in $W^u(x)$.
A very simple, but crucial, observation is that the
semi-flow $\varphi_s$ deforms the \emph{ascending disk}
$W^s_\eps(x):=W^s(x)\cap \Lambda^{c+\eps} M=W^s(x)\cap N_x$
to $x$, as $s\to\infty$. Clearly this fails on other parts
of $N_x$. Note that $W^s_\eps(x)$ is a $C^1$ graph over its
tangent space denoted by, say $X^+$. The idea is to
\emph{foliate all of $N_x$ by copies of $W^s_\eps(x)$, more precisely $C^1$ graphs over $X^+$, then extend $\varphi_s$ artificially to all of $N_x$ using the graph maps};
see~(\ref{eq:induced-semi-flow}) and Figure~\ref{fig:fig-induced-flow-JP}.

\begin{figure}
  \centering
  \includegraphics{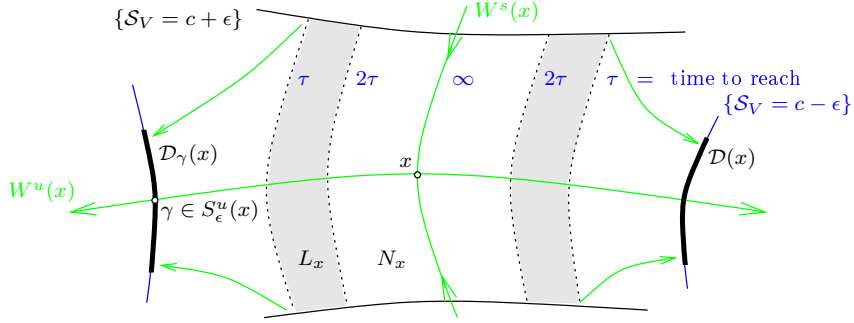}
  \caption{Conley pair $(N_x,L_x)$ foliated by
           equal time disks $(\varphi_T)^{-1}\Dd_\gamma(x)$}
  \label{fig:fig-N-foliated-JP}
\end{figure}

To see the foliation assign to each point of $N_x$ the time
$T$ at which it hits the level surface $\{\Ss_V=c-\eps\}$; see
Figure~\ref{fig:fig-N-foliated-JP}.
This suggests that $N_x$ is foliated
by (pieces of) the equal time hypersurfaces
$(\varphi_T)^{-1}\{\Ss_V=c-\eps\}$ for $T\in(\tau,\infty)$.
For $T=\infty$ one obtains the codimension $k$
ascending disk $W^s_\eps(x)$.
Of course, the leaves of a foliation need to be of the
same codimension: Consider the tubular neighborhood
$\Dd(x)\to S^u_\eps(x)$ associated to the (sufficiently small)
radius $a$ normal bundle of the descending sphere
$S^u_\eps(x):=W^u(x)\cap\{\Ss_V=c-\eps\}$
in the Hilbert manifold $\{\Ss_V=c-\eps\}$.
Each fiber $\Dd_\gamma(x)$ is a codimension $k$ disk.

\section{Backward \boldmath$\lambda$-Lemma and stable foliations}\label{sec:backward-lambda}

Fix $x\in\Crit_k$ and set $c:=\Ss_V(x)$.
Since $N_x=N_x^{\eps,\tau}$ fits into any neighborhood
of $x$ for $\eps>0$ small and $\tau>0$ large
we use local coordinates about $x\in\Lambda M$.

\subsubsection*{Local coordinates about \boldmath$x\in\Lambda M$}
The nonlinear part of the heat equation~(\ref{eq:heat})
determines a closed radius $\rho_0$ ball $\Bb_{\rho_0}$
about $0\in X$ such that the following is true.
Paths $s\mapsto u(s)$ in $\Lambda M$ near
$x$ and $s\mapsto\xi(s)$ in $\Bb_{\rho_0}$ uniquely
correspond to each other via the identity
$u(s)=\exp_x\xi(s)$ pointwise for every $t\in S^1$.
In the new coordinates $\xi$ the Cauchy problem associated
to~(\ref{eq:heat}) turns into the equivalent Cauchy problem
\begin{equation}\label{eq:cauchy-local}
     \zeta^\prime(s)+A\zeta(s)
     =f(\zeta(s)),\qquad
     \zeta(0)=z\in\Bb_{\rho_0},
\end{equation}
for maps $\zeta:[0,T]\to\Bb_{\rho_0}\subset X$. Here $A=A_x$
is the Jacobi operator associated to the (perturbed) closed geodesic
$x$. The semi-flow $\varphi$ turns into a local
semi-flow $\phi$ on $\Bb_{\rho_0}\subset X$. The nondegenerate
critical point $x$ corresponds to the hyperbolic fixed point
$0$ of $\phi$. Furthermore, there is the orthogonal splitting
$$
     X
     :=T_x\Lambda M
     \simeq T_xW^u(x)\oplus T_xW^s(x)
     =: X^-\oplus X^+.
$$
Here $X^-$ is of finite dimension $k=\IND_V(x)$
and consists of smooth loops along $x$.
By $\pi_\pm:X\to X^\pm$ we denote the associated orthogonal
projections. For coordinate representatives
of global objects we shall use the global notation omitting $x$,
for example $W^u(x)$ becomes $W^u$. By $\Ss$ we
denote the representative of $\Ss_V$.
Via a (standard) change of coordinates one achieves
that locally near zero $W^u$ is contained in $X^-$.
By $\Bb^+_R$ we denote the radius $R$
ball about $0\in X^+$. The
\emph{spectral gap} $d>0$ is the distance between $0$ and
the spectrum of $A_x$.

\begin{figure}
  \centering
  \includegraphics{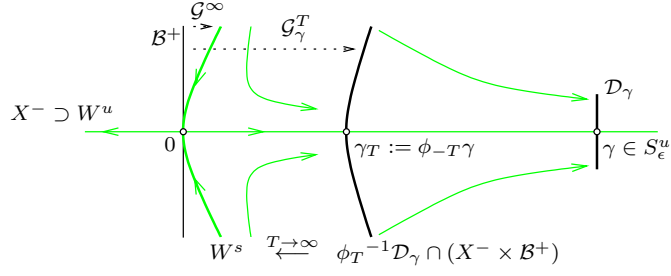}
  \caption{Backward $\lambda$-Lemma}
  \label{fig:fig-lambda-lemma}
\end{figure}

\begin{theorem}[Backward $\lambda$-Lemma, \cite{Joa-LAMBDA}]
\label{thm:backward-lambda-lemma}
Pick $\mu\in(0,d)$
and a hypersurface $\Dd\subset\Bb_{\rho_0}$ of the form
$S^u_\eps\times\Bb^+_a$. Then the following is true
(see Figure~\ref{fig:fig-lambda-lemma}).
There is a ball $\Bb^+$ about $0\in X^+$, a constant
$T_0>0$, and a Lipschitz continuous map
\begin{equation*}
\begin{split}
     \Gg:(T_0,\infty)\times S^u_\eps\times\Bb^+
    &\to W^u\times\Bb^+\subset\Bb_{\rho_0}
   \\
     (T,\gamma,z_+)
    &\mapsto 
     \left(G^T_\gamma(z_+),z_+\right)
     =:\Gg^T_\gamma(z_+)
\end{split}
\end{equation*}
of class $C^1$. Each map $\Gg^T_\gamma:\Bb^+\to X$
is bi-Lipschitz, a diffeomorphism onto its image, and
$\Gg^T_\gamma(0)=\phi_{-T}\gamma=:\gamma_T$.
The graph of $G^T_\gamma$ consists of those $z\in \Bb_{\rho_0}$
which satisfy $\pi_+ z\in \Bb^+$ and reach the fiber 
$\Dd_\gamma=\{\gamma\}\times\Bb^+_a$
at time $T$, that is
$$
     \Gg^T_\gamma(\Bb^+)
     =(\phi_T)^{-1}
     \Dd_\gamma\cap
     \left( X^-\times\Bb^+\right).
$$
Furthermore, the graph map $\Gg^T_\gamma$ converges
uniformly, as $T\to\infty$, to the stable manifold
graph map $\Gg^\infty$.
More precisely, the estimates
\begin{gather*}
     \Norm{\Gg^T_\gamma(z_+)-\Gg^\infty(z_+)}_{W^{1,4}}
     \le e^{-T\frac{\mu}{16}}
     ,\qquad
     \Norm{d\Gg^T_\gamma(z_+)v}_2
     \le 2 \Norm{v}_2,
     \\
     \Norm{d\Gg^T_\gamma(z_+)v-d\Gg^\infty(z_+)v}_2
     \le e^{-T\frac{\mu}{16}}\Norm{v}_2
\end{gather*}
hold true for all $T>T_0$, $\gamma\in S^u_\eps$,
$z_+\in\Bb^+$, and $v$ in the $L^2$ closure of $X^+$.
\end{theorem}
Theorem~\ref{thm:backward-lambda-lemma} is based on the
observation that the Cauchy problem for a heat flow line
$\xi:[0,T]\to X$ with $\xi(0)=z$
is equivalent to a \emph{mixed Cauchy problem} with data
$(T,\gamma,z_+)$. Namely, there is a unique heat flow line
$\xi:[0,T]\to X$ with $\pi_+\xi(0)=z_+$ and
$\pi_-\xi(T)=\gamma$. 

That the  ($k$-dimensional) unstable manifolds
carry backward time information is evident
from their definition.
In contrast,
Theorem~\ref{thm:backward-lambda-lemma} provides
backward time information on \emph{open} sets.

\subsubsection*{Stable foliation of Conley set}

Theorem~\ref{thm:backward-lambda-lemma}
foliates neighborhoods of $x$ by (globally meaningless)
codimension $k$ disks. The next result provides global
information in various directions.
By definition the \emph{descending disk} $W^u_\eps(x)$ is
given by $W^u(x)\cap\{\Ss_V>c-\eps\}$.
\begin{theorem}[\cite{Joa-CONLEY}]\label{thm:inv-fol}
Given $\mu\in(0,d)$ there are constants $\eps_1,\tau_1,a>0$
such that the following is true. Assume $\tau>\tau_1$ and
$\eps\in(0,\eps_1)$ and consider the radius $a$ tubular
neighborhood $\Dd(x)\to S^u_\eps(x)$ defined in the paragraph
preceding section~\ref{sec:backward-lambda}.
\begin{enumerate}
\item[a)]
  The Conley set $N_x=N_x^{\eps,\tau}$ carries the structure
  of a codimension $k$ foliation whose leaves are
  parametrized by the disk
  $\varphi_{-\tau} W^u_\eps(x)$.
  The leaf over $x$ is the
  ascending disk $W^s_\eps(x)$ and the other leaves
  are given by the disks
  \begin{equation*}
     N_x(\gamma_T)
     =\left\{(\varphi_T)^{-1}\Dd_\gamma (x)\cap
     \{\Ss<c+\eps\}\right\}_{\gamma_T}
     ,\qquad
     \gamma_T:=\varphi_{-T}\gamma,
  \end{equation*}
  whenever $T>\tau$ and $\gamma\in S^u_\eps(x)$.
\item[b)]
  Leaves and semi-flow are compatible in the sense that
  \begin{equation*}
     z\in N_x(\gamma_T)\quad\Rightarrow\quad
     \varphi_\sigma z\in N_x(\varphi_\sigma\gamma_T)
    ,\quad
     \forall \sigma\in[0,T-\tau).
  \end{equation*}
\item[c)]
  The leaves converge uniformly to the ascending disk
  in the sense that
  \begin{equation*}
     \dist_{W^{1,2}}\left(N_x(\gamma_T),W^s_\eps(x)\right)
     \le e^{-T\frac{\mu}{16}}
  \end{equation*}
  for all $T>\tau$ and $\gamma\in S^u_\eps(x)$.
  Furthermore, if $U$ is a $\delta$-neighborhood
  of $W^s_\eps(x)$ in $\Lambda M$, then
  $N_x^{\eps,\tau_*}\subset U$ for some constant $\tau_*$.
\item[d)]
  Assume $U$ is an open neighborhood of $x$ in
  $\Lambda M$. Then there are constants $\eps_*$ and
  $\tau_*$ such that $N_x^{\eps_*,\tau_*}\subset U$.
\end{enumerate}
\end{theorem}

\section{Strong deformation retract}

Pick $x\in\Crit_k$. It remains to prove~(\ref{eq:claim}).
If $k=0$, then $L_x=\emptyset$ and $W^u(x)=\{x\}$ is a strong
deformation retract of $W^s_\eps(x)=N_x$. The retraction is
provided by the semi-flow $\varphi_s$ and we are done.
Assume $k>0$. Consider the local
setup of section~\ref{sec:backward-lambda} and
denote the representative of $N_x$ by $N$;
similarly for other quantities.
Fix $\rho_0>0$ so small that the only
critical point in $\Bb_{\rho_0}$ is $0$.

\begin{definition}\label{def:induced-semi-flow}
By Theorem~\ref{thm:inv-fol}
each $z\in N$ lies on a leaf $N(\gamma_T)$ for
some time $T>0$ and some point $\gamma$ in the
descending disk $S^u_\eps$ where
$\gamma_T:=\phi_{-T}\gamma$.
The continuous leaf preserving map
$\theta:[0,\infty)\times N\to N$ defined by
\begin{equation}\label{eq:induced-semi-flow}
     \theta_s z
     :=\Gg_\gamma^T
     \pi_+ \phi_s
     \Gg^\infty \pi_+ z
\end{equation}
is called the
\emph{induced semi-flow on $N$};
see Figure~\ref{fig:fig-induced-flow-JP}.
It is of class $C^1$ on $(0,\infty)\times N$.
\end{definition}
\begin{figure}[b]
  \centering
  \includegraphics{fig-induced-flow-JP.eps}
  \caption{The induced flow $\theta_s$ on $N$}
  \label{fig:fig-induced-flow-JP}
\end{figure}

That $\theta_s$ preserves the central leaf $N(0)=W^s_\eps$
is due to the downward $L^2$-gradient
nature of the heat equation. The proof for a general
leaf $N(\gamma_T)$ turns out to be surprisingly complex
although the idea is once more simple:
Show that the map $s\mapsto\Ss(\theta_s z)$ strictly
decreases whenever $z$ lies in the (topological) boundary
of a leaf. This implies preservation of leaves as follows.
Firstly, note that $\theta$ is actually defined on a
neighborhood of $N(\gamma_T)$ in $\Gg^T_\gamma(\Bb^+)$.
Secondly, the (topological) boundary of a leaf
lies on action level $c+\eps$ whereas the leaf itself
lies strictly below that level. Thus the induced
semi-flow points inside along the boundary of each
leaf---which is a disk by Theorem~\ref{thm:inv-fol}. So
$\theta_s$ preserves leaves, thus $N$ and $L$ by
Theorem~\ref{thm:inv-fol}.
Moreover, it continuously deforms both topological spaces
to their respective part in the unstable manifold and
this concludes the proof of~(\ref{eq:claim}).
Therefore $\Ff$ defined by~(\ref{eq:F_k})
is indeed a Morse filtration for $\Lambda^a M$ and
by Remark~\ref{rem:Morse-filtration} this establishes the
desired natural isomorphism~(\ref{eq:natural-isomorphism}).

\begin{figure}
  \centering
  \includegraphics{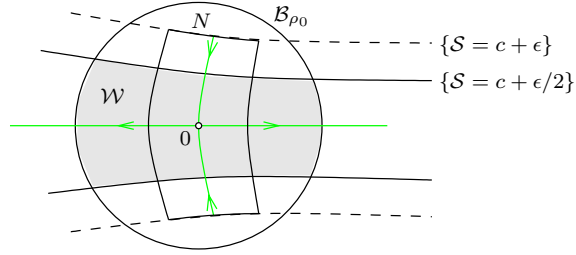}
  \caption{The neighborhood $\Ww$ of $0$
           used to define $\alpha>0$}
  \label{fig:fig-neighborhood-Ww-JP}
\end{figure}

It remains to show that $\frac{d}{ds}\Ss(\theta_s z)<0$
whenever $z$ lies in the (topological) boundary of a leaf.
Note that $\grad\Ss$ is defined on loops whose regularity
is at least $W^{2,2}$. Consider the neighborhood
$\Ww:=\Bb_{\rho_0}\cap\{\Ss\le c+\eps/2\}$ of $0$ illustrated
by Figure~\ref{fig:fig-neighborhood-Ww-JP}.
By Palais-Smale the constant defined by
\begin{equation*}
     \alpha:=
     \inf_{z\in\left(\Bb_{\rho_0}\cap W^{2,2}\right)\setminus \Ww}
     \Norm{\grad\Ss(z)}_2
     >0
\end{equation*}
is strictly positive.
A rather technical argument, see~\cite{Joa-CONLEY}, involving
a long calculation which uses heavily the estimates provided by
Theorem~\ref{thm:backward-lambda-lemma} shows that for all
$\eps>0$ small and $\tau>0$ large the following is true.
If $T>\tau$ and $\gamma\in S^u_\eps$, then
\begin{equation*}
\begin{split}
     \frac{d}{ds}\Ss(\theta_s z)
    &=d\Ss|_{\theta_s z} \,
     d\Gg_\gamma^T|_{z_+(s)} \,
     \pi_+  
     \tfrac{d}{ds}
     \left(\phi_s\Gg^\infty\pi_+ z\right)
   \\
    &=-\left\langle
     \grad\Ss|_{\theta_s z},
     d\Gg_\gamma^T|_{z_+(s)}
     \pi_+
     \grad\Ss|_{\phi_s q}
     \right\rangle_{L^2}
   \\
    &\le -\tfrac{1}{4}\alpha^2
\end{split}
\end{equation*}
for all $z\in\p N(\gamma_T)$ and $s>0$ small.
It is precisely this calculation where we need 
convergence in $W^{1,4}$ and the
extension to $L^2$ of the linearized graph map
$d\Gg_\gamma^T(z_+)$ in
Theorem~\ref{thm:backward-lambda-lemma}.
(The nonlinear part $f$ of~(\ref{eq:heat})
maps $W^{1,4}$ to $L^2$.)

\vspace{.2cm}
\noindent
{\bf Acknowledgements.} The author would like to thank Alberto Abbondandolo, Kai Cieliebak, and Klaus Mohnke for most useful discussions and comments.


\end{document}